\documentclass[12pt]{article}
\usepackage[numbers]{natbib}
\usepackage{amssymb}
\begin{document}

\title{\bf Dynamical real numbers and living systems
}

\author{
Dhurjati Prasad Datta\thanks{email:dp${_-}$datta@yahoo.com} \\
Department of Mathematics,
North Bengal University,\\
 P.O. North Bengal University, 
Darjeeling, India, Pin: 734430} 
\date{}
\maketitle

\baselineskip = 20pt 

\begin{abstract}
Recently uncovered second derivative discontinuous solutions of the simplest linear 
ordinary differential equation define not only an nonstandard extension of the framework 
of the ordinary calculus, but also provide a dynamical representation of the ordinary 
real number system. Every real number can be visualized as a living cell -like 
structure,  endowed with a definite evolutionary arrow. We discuss the relevance of this 
extended calculus in the study of living systems. We also present an intelligent version 
of the Newton's first law of motion.   
\end{abstract}

\begin{center}
\em{Chaos, Solitons, and Fractals, vol 20, issue 4, 705-712, (2004)}
\end{center}

\newpage

\section{Introduction}
Let $t$  be a real variable. Then a key assumption in the formalism of mathematical 
analysis (calculus), and hence in classical and quantum dynamics, is that $t$ changes 
only by (linear) translations. The variable $t$ assume values from the set of real 
numbers, imagined as fixed (static), ideal mathematical objects. The fundamental concept 
of limit, as treated in the standard $\epsilon,\delta$-definition,  disregards the 
inherent dynamical aspect of the concept, by reducing it to a relative statement between 
`sufficiently small' real numbers. The growth of the modern (standard)  analysis, thus 
not only hides the deep interrelationship between analysis and dynamics, but also 
rejects infinitesimally small (infinitely large) numbers. The nonstandard analysis of 
Robinson ~\cite{robin} , although puts the theory of infinitesimally small numbers into 
a sound mathematical basis with interesting applications in the theory of fractals and 
stochastic analysis ~\cite{nottale}, the true dynamic features of the infinitesimals 
have not been appreciated. The recently discovered second derivative discontinuous scale 
free solutions of the linear ordinary differential equations (of ordinary 
calculus)~\cite{dp1,dp2} now provides, in one hand, a new direct proof of the existence 
of infinitely small numbers in the set of real numbers, and on the other, shows that 
such numbers must be `intrinsically dynamic' in the sense that they are in a state of 
spontaneous random fluctuations, which, in turn, induces a sense of evolution to every 
ordinary real number. The qualifier  `intrinsic' means that the said evolution is not 
(actually it can not be) due to any externally applied force, as it is necessary in the 
treatment of ordinary dynamics, but instead an intrinsic property of the real number 
system itself. In this paper, we point out a few more interesting features of the 
dynamical number system and argue that the extended `calculus of dynamical numbers' 
should be the ideal framework for the study of living (intelligent) systems. We also 
give an intelligent version of the Newton's first law of motion, which should be 
considered as the `first law of motion' for living systems. Let us recall here that the 
formalism of calculus (and dynamics) was originally invented to study the motion of 
material bodies, the planetary system, for instance. The modern studies of living 
systems (life sciences) tend in fact to reduce `life' essentially to biochemical 
(physical) processes, although there have been occasional suggestions for the need of an 
intelligent principle, rather than simply material (physical) principles as emerged from  
the studies of quantum field theories and nonlinear sciences, to understand the complex 
dynamics of a living system~\cite{roger}. The main theme of the present paper is 
therefore to highlight the potentially interesting dynamical features of the new scale 
free calculus in the context of a living intelligent system. 
\par The singular role of fractals in the studies of nonlinear systems is now well 
appreciated. Recently, there have been a lot of studies devoted to applications of 
fractal geometry, deterministic chaos and stochastic dynamics in classical and quantum 
dynamics~\cite{ord, nottale, enash1, enash2, sidh}, as well as in biological systems 
~\cite{west}. However, the inherent dynamic nature of any fractal object, though 
implicitly present in their definitions ( as the limit set of an infinite iterative 
process, inducing a temporal sense), is yet to be enunciated clearly. Our results 
indicate that any dynamical system, being inherently scale free and directed, must tend 
toward a random fractal like evolutionary process, with approximate (statistical) 
self-similarity over multiple random scales. Conversely, the emergence of fractal 
objects must be viewed in an evolutionary setting, provided naturally by the scale free 
extension of calculus. Our work also tells that points in the real number line, and 
hence in the space-time manifold of any physical theory, have a Cantorian structure, but 
this need not necessarily be linked only with the Planck scale physics, as in the 
$\mathcal E^{\infty}$ theory ~\cite{enash1}, but may arise more universally and 
naturally to any truly dynamical theory, for instance, in the context of living systems.          
The possibility of random fractal topology, as an inherent manifestation of stochastic 
dynamics near any fixed point of the underlying quantum field theory, has recently been 
pointed out ~\cite{egold,enash3,enash4}.

\section{Scale free calculus}
\subsection{Generalized Solutions}
\par Let us begin by reviewing the important results of the new scale free calculus 
~\cite{dp1, dp2} . In    ordinary (linear) calculus, one begins with an independent real 
variable $t$ (say) and then introduces the concepts of limit, continuity, derivatives 
and so on for a function of $t$. For definiteness, the variable $t$ is said to define 
the O(1) scale of ordinary (macroscopic/ coarse grained) real numbers. The scale free 
extension of calculus, however, is defined on the basis of scale free, second derivative 
discontinuous solutions of the linear differential equation 

\begin{equation}\label{si}
{\frac{{\rm d}T}{{\rm d}t}}=T
\end{equation}

\noindent which is given by 

\begin{equation}\label{gen}
\ln T(t)=t +r k \phi(t_1),\,\phi=t_1 \tau(t_1^{-1})
\end{equation}

\noindent where $\tau(t_1)=T(\ln t_1),\,t_1=k_0t,\,k_0>0$ being a solution of 
eq~(\ref{si}), but in the log scale, denotes a self~-similar replica of $T(t)$ in the 
smaller logarithmic variable $\ln t_1^{-1}$ and $k$ and $k_0$ are arbitrary scaling 
parameters. Here $r$ is a random variable assuming values from $\{+1,-1\}$ with equal 
probability, say. The solution is scale free, since both $\tau$ and $\phi$ are solutions 
of scale free equations 
\begin{equation}\label{ss}
t{\frac{{\rm d}\tau}{{\rm d}t}}=\tau
\end{equation}
\noindent and 
\begin{equation}\label{sc}
t{\frac{{\rm d}\phi}{{\rm d}t}}=0
\end{equation}

\noindent respectively. Clearly, $T(t)$ satisfies eq(\ref{si}) if and only if $\phi$ 
satisfies eq(\ref{sc}), which means, in turn, that $\tau$ must satisfy eq(\ref{ss}).   
The scaling parameter $k_0$ is representative of this scale free nature. Utilizing the 
scaling law $\tau(kt)=k\tau(t)$, and the scale free property, one can always re-scale 
$\tau$ suitably to make $k$ arbitrarily small. Further, the generalized solution is 
singular at $t=0$ and intrinsically nonlocal since the implicit definition of the 
generalized function $T$ involves distinct points  $t$ and $\ln t_1^{-1}$, for instance, 
over two different scales.  These two properties viz., {\em the singularity and 
nonlocality, tell that $T(t)$ is a new, nontrivial solution of eq(\ref{si})}. (The well 
known standard solution of eq~(\ref{si}) $\ln T_s(t)=t +c,\, c$ being an ordinary 
constant, is local (i.e., defined at a single point) and well defined for all $t$ in the 
real number set $R$. Moreover, replacing $\tau$ in eq(\ref{gen}) by the standard 
solution $\tau=c t$ (say), one reproduces the standard solution $T_s$, but in  $R-\{0 
\}$. Consequently, the standard solution belongs to this more general class of 
solutions, when the initial condition at $t=0$ is defined only approximately. In other 
words, the standard solution is only an approximate solution, in the context of this 
more general class of solutions.   Further, the framework of ordinary calculus does not 
allow randomness in the context of eq(\ref{si}).)  The nontrivial scale free solution 
$\phi$ of eq~(\ref{sc}) is called a slowly varying generalized constant, since $\phi$ 
remains constant for any finite $t$, but may experience a slow variation on the 
infinitesimal scale $\ln t_1^{-1}$, for an arbitrarily large $t$.
\par The nontrivial part of the generalized solution $\ln T$ corresponds to {\em 
infinitesimally small} random scales. In fact, even as the continuous (ordinary) real 
variable $t$ exhausts all arbitrarily small ordinary real numbers in the vicinity of 
$t=0$, $\ln T(0)=\lim_{t\rightarrow 0}k\phi(t) $ remains arbitrarily small but nonzero, 
for an arbitrarily small but nonzero generalized constant $k$. Note that $k$ can not 
vary faster than $\phi$, by definition. Further, it is useful to model the arbitrariness 
of the scale factor $k$ by a continuous random variable $K$ which can assume values $k$ 
with a uniform probability density function, for instance. However, as shown in 
Ref\cite{dp1}, a Gamma-like distribution seems to be an appropriate distribution in the 
context of the above scale free solutions.
\par We note, further, that the singularity of $\phi$ is not specific to $t=0$ only. In 
fact, by the translation symmetry of eq(\ref{si}), $\phi(t^\prime)=\phi(t),
\,t^\prime=t-t_0$, so that (i) $\phi(t)(=1+k\phi(t_1))$ is actually singular at every 
point $t_0$ in the ordinary real axis and (ii) is essentially a universal random 
function in the neighbourhood of 1 ($k\phi$ is infinitesimal in the neighbourhood of 0), 
where both $t$ and $t_1\,(t\neq t_1)$ are random variables defined near 1. Thus, 
$\lim_{t\rightarrow t_0}(\ln T(t)-t_0)=\ln T(0)(\neq 0)$. In the ordinary calculus, this 
limit must necessarily equal to 0. 

{\em Definition 1}:  As the limit $t\rightarrow 0$ exhausts all the coarse grained real 
numbers, there exists a continuum of random, infinitesimally small {\em microscopic} 
positive numbers (scales)  given by  $\eta= k\phi(t_1),\,k>0  $. Clearly, $\eta\neq 0$ 
is less than any arbitrarily small, nonzero, coarse grained, positive real number. 

{\em Remark 1:} The standard solution of eq~(\ref{si}) $\ln T_s(t)=t, T(1)=e$ defines a 
1-1 (identity) mapping of the real number set $R$ onto itself. The generalized solution 
$\ln T(t)$ now defines an extension of the identity mapping from  $R$ to the nonstandard 
real set {\bf R}. Clearly, {\bf R}=$R$, $T$ being an exact solution of eq~(\ref{si}). 
Accordingly, every ordinary real number $t$ is identified (extended) with (to) the set 
(actually an equivalence class under the relation `infinitely close neighbour') of scale 
free (fat) real numbers denoted ${\bf t}=t + r N_0 $, where $N_0$ is the nontrivial 
neighbourhood of positive infinitesimals $N_0=\{\eta \}$. The macroscopic zero 0, for 
instance, corresponds to the set of scale free (first order) infinitesimals {\bf 0}=$\{r 
\eta\} (<{\bf 0>}=0$). Here, $< >$ denotes the expectation value of the random variable 
{\bf t}. Determining the exact value of a real number is thus theoretically impossible, 
because of the identity {\bf R}=$R$. Further, the additive zero of ${\bf R}$ now 
corresponds to the set of higher order infinitesimals, viz.,  $0=\{\eta^2,\, \eta^3, 
\ldots\}$. 

{\em Remark 2:} The coarse grained variable $t$ undergoes changes by linear 
translations. The intrinsic randomness of the infinitesimals now tells that the 
infinitely small numbers are in a state of continual fluctuations, analogous to a system 
of ideal gas in thermodynamic equilibrium, say. The fluctuations in the neighbourhood of 
the normalized class ${\bf 1}\equiv{\bf t}/t$ could be modeled by the inversion rule 
$t_{-}\rightarrow t_{-}^{-\alpha}=t_{+}^{\alpha}$ where $\alpha > 0$ and $t_{\pm}=1 \pm 
\eta,\, \eta$ being an (first order) infinitesimal so that O($\eta^2)=0$, by definition. 
The points in the infinitesimal neighbourhood of {\bf 1} can therefore change by 
discrete jumps (inversions), the length of  jumps being random because of the 
arbitrariness of $\alpha$. That inversion transformations should be operative here is 
actually traced to the non-locality of the generalized solutions (for a more direct 
proof see below). These small scale fluctuations are, however, un-noticeable  at the 
level of the macroscopic scale, so that $t=<{\bf t}>$ behaves as a smooth linear 
variable. We also note that $<\phi(t)>=<\phi(t^{-1})>=<{\bf 1}>=1$.

\par To prove the second derivative discontinuity of the generalized solution and to 
establish its directed and irreversible character, let us now present a simpler  
solution of eq(\ref{ss}) displaying the possibility of inversion in an explicit manner.     

\subsection{Second derivative discontinuity and Irreversibility}

Let $t_{\pm}=1 \pm \eta,\, 0<\eta<< 1$, so that O($\eta^2)=0$ and $t_-\rightarrow 
t_-^{-1}=t_+$.
Let also that $\tau (t_{\pm})=\tau_{\pm}$. Then the function $\tau_f(t)$
defined by  
\begin{equation}\label{ns}
\tau_f(t)=
\cases{
\tau_- & {\rm if} $t\lessapprox 1$\cr
\tau_+ & {\rm if} $t\gtrapprox 1$},\,\,
\tau_{-}(t_-)=1/t_+,\, \tau_+(t_+)=t_+
\end{equation}

\noindent constitutes a new nontrivial solution of eq(\ref{ss}), in the 
neighbourhood of $t=1$, with the initial condition $\tau_f(1)=1$. Indeed, it is easy to 
verify, using ${\rm d}t_-=-{\rm d}t_+$ (and/or $t_-^{-1}{\rm d}t_-=-t_+^{-1}{\rm d}t_+$, 
since $\ln t_-=-\ln t_+$),  that  

\begin{equation}\label{fd}
{\frac{{\rm d}\tau_{-}}{{\rm d} t_-}}=-{\frac{{\rm d}\tau_{-}}{{\rm d} t_+}}=1/t_+^2,\, 
{\frac{{\rm d}\tau_{+}}{{\rm d} t_+}}=1
\end{equation}

but 

\begin{equation}\label{sd}
{\frac{{\rm d}^2\tau_{-}}{{\rm d} t_-^2}}=2/t_+^3,\, {\frac{{\rm d}^2\tau_{+}}{{\rm d} 
t_+^2}}=0
\end{equation}

\noindent which tell that the function $\tau_f(t)$ defined by eq(\ref{ns}) has 
continuous first derivative at and in the vicinity of $t=1$, but its second derivative 
has a discontinuity at $t=1$. One notes, moreover, that eq(\ref{fd}) can be written as   
\begin{equation}\label{tsl}
t_-{\frac{{\rm d}\tau_{-}}{{\rm d} t_-}}=\tau_-,\, t_+{\frac{{\rm d}\tau_{+}}{{\rm d} 
t_+}}=\tau_+
\end{equation}

\noindent so that $\tau_f(t)$ is a solution. Because of its second derivative 
discontinuity, it is indeed a {\em new} solution, distinct from the infinitely 
differentiable standard solution. 

\par A fundamental property of the new solution eq(\ref{ns}), which is {\em exact} in 
the extended {\bf R} with nontrivial infinitesimals, is that it breaks the reflection 
(time reversal/parity) symmetry of the underlying equation (\ref{ss}). (The solution is 
approximate (O($\eta^2$)) in the ordinary real number set $R$. The origin of the second 
derivative discontinuity is then traced to the approximation introduced by the 
substitution $t_-^{-1}=t_+$. In ${\bf R}$, however, this approximation is an intrinsic 
property of the infinitely small numbers, which can not be improved to an unlimited 
accuracy. In the Appendix, we present an exact (to all orders of $\eta$) random solution 
of eq(\ref{ss}).) Let $P:\, Pt_\pm=t_\mp$ denote the parity transformation  close to 
$t=1$. Then the parity reversed function $\tau_f^P: \,\tau_-^P=t_-,\, \tau_+^P=1/t_-$ is 
also a solution of eq(\ref{ss}), but, clearly it is different from the original 
solution $\tau_f$. We note, in particular, that $\tau_f^P$ diverges as $\eta\rightarrow 
1^-$, in contrast to $\tau_f$ which is well-behaved, as it should be. Consequently, one 
expects that the solution $\tau_f$ is more likely to materialize, compared to the parity 
reversed solution. As a result, it is more natural for the variable $t$ to change from 
$t_-$ to $t_+$ (by an instantaneous jump) via $t_-\rightarrow t_-^{-1}=t_+$, since the 
solution (\ref{ns}) admits the transition $\tau_-\rightarrow \tau_+=1/\tau_-$, and then 
is forced to increase with $\eta\rightarrow 1^-$ until a second generation inversion in 
the small scale variable $\eta$, rather than to change alternatively. Consequently, the 
solution eq(\ref{ns}) tells in a succinct way that the function $\tau_f(t)$ must always 
be monotonically increasing and hence is directed irreversibly  in the forward 
direction.    
       
\par We note that the problem of time asymmetry (arrow of time) in Physics is that all 
the fundamental equations of motion are time reversal symmetric. Consequently, one 
expects that the corresponding solutions must also respect the symmetry of the 
underlying equation and so must be time reversal symmetric. The present solution is the 
{\em first explicit} example which fails to respect the symmetry of the underlying 
equation, analogous to the spontaneous (vacuum) symmetry breaking encountered in the 
quantum field theory. As advocated many a times by Prigogine~\cite{prigo} (and also 
implicit in El Naschie's~\cite{enash2} work), any time reflection symmetric system must 
essentially be interpreted as devoid of an intrinsic time sense. In this sense a state 
in thermal equilibrium is timeless since it lacks an intrinsic tendency of evolution 
toward a final state. Accordingly, a true time variable should be one which breaks this 
reflection symmetry spontaneously. Our solution naturally stands out as a {\em genuine 
candidate for time}.    

{\em Remark 3:} The solution (\ref{ns}) shows explicity that inversion $t_\mp=1/t_\pm$ 
is as natural a mode of `increment' (change) for a real variable, as it is usually 
assigned to translation. Recall that the standard solution of eq(\ref{ss}) under pure 
translation (i.e., when $t$ is allowed to change only by translation) is 
\begin{equation}\label{stand}
 \tau (t_-)=t_-,\, \tau(t_+)=t_+
 \end{equation}
 
 \noindent We note here also that another exact solution of eq(\ref{ss}) (under 
inversion) is 
\begin{equation}\label{fluc}
\tau_-=1/t_+,\, \tau_+=1/t_-
\end{equation}

\noindent Comparing the above equation with eq(\ref{ns}), we see that in eq(\ref{ns}) 
only `one sided inversion' is realized, whereas eq(\ref{fluc}) accommodates inversions 
in both the ways. Moreover, solution (\ref{fluc}) is an effect of pure inversion only, 
when eq(\ref{ns}) is generated by accommodating both the translation and inversion modes 
of change for $t$. One also verifies that eq(\ref{fluc}) is infinitely differentiable, 
as in eq(\ref{stand}), and is time reflection symmetric. The solution (\ref{fluc}) thus 
has the natural status of {\em `equlibrium' infinitesimal fluctuations} near $t=1$. We 
recall that the equilibrium fluctuations in thermodynamics is time symmetric (c.f., 
remark 2).

{\em Remark 4:} A more general class of time asymmetric solution of eq(\ref{ss})  
 
\begin{equation}\label{as}
\tau^\alpha_f(t)=
\cases{
\tau_- & {\rm if} $t\lessapprox 1$\cr
\tau_+ & {\rm if} $t\gtrapprox 1$},\,\,
\tau_{-}(t_-)=t_+^{-\alpha},\, \tau_+(t_+)=t_+, \, \alpha > 0
\end{equation}

\noindent when $t_-$ relates to $t_+^{\prime\alpha}=1+\alpha\eta^\prime >1$ via 
$t_-\rightarrow  t_-^{-\prime\alpha}=t_+^{\prime\alpha}$ so that ${\rm d}\eta={\rm 
d}\eta^\prime$.\footnote{ The case $\alpha < 0$ is also interesting. Two points 
$t_1=1+\eta_1$ and $t_2=1+\eta$, on the same side of 1, may be related by an 
infinitesimal scaling $t_1=t_2^{\beta},\,\eta_1=\beta\eta,\,\alpha=
-\beta <0$. The corresponding solution is $\tau(t_1)=t_2^{\beta},\, \tau(t_2)=t_2$.} 
Clearly, the solution (\ref{as}) reduces to eq(\ref{ns}) for $\alpha=1$. 
We note that a power law reduces to a 
scaling (dilaton) transformation for small enough $\eta: \,\eta_1=\alpha 
\eta$ near $t=1$, satisfying $0<\alpha\eta<<1$.  
One can also verify that the solution (\ref{as}) is only 
first order differentiable, the second derivative having a discontinuity 
at $t=1$. It also follows that a real variable 
$t$ can not only change, in the vicinity of $t=1$, by a pure inversion 
$\alpha=1$, but also undergo changes by discrete random jumps, 
length of a jump being determined by the arbitrary parameter $\alpha$.
Thus a real variable $t(\lesssim 1)$, (viz. $t_-=1-\eta$), for instance,  can 
now connect instantaneously to a spectrum of values $t_{\alpha+}=1+\alpha \eta$, 
which can naturally be assumed to be distributed following a probability law.
In other words, the existence of this class of solutions puts our original assumption of 
fluctuating infinitesimals in the neighbourhood of a real number ~\cite{dp1} in a farmer 
footing.  

\par The above class of time asymmetric solutions now motivates the 

{\em Definition 2:} Let $T(t)=t\phi(t_1)=t(1+k t_1^{-1}\tau_f(t_1))$. The variable $T$ 
is intrinsically time-like if and only if $\tau_f$ belongs to the class of solutions 
given by eq(\ref{ns}) (equivalently, eq(\ref{as}).

{\em Remark 5:} In case $\tau_f$ belongs either to the solution (\ref{stand}) or 
(\ref{fluc}) the corresponding $T$ would be reflection symmetric (as well as infinitely 
differentiable) and can not be considered as a genuine time variable. We note that the 
contemporary dynamical theories make use of a time variable belonging to the solution 
(\ref{stand}) only. We note in this context that the conventional notations $t$ and $x$ 
for time and space variables can not distinguish the two types of variables in an 
intrinsic way, since both the variables are reflection symmetric solution (\ref{stand}) 
of eq(\ref{ss}).

\par As a consequence of the above remarks, one can restate the definition 1 thus

{\em Definition 3:} A variable $t$ is macroscopic if it belongs to the trivial 
reflection symmetric solution (\ref{stand}). A variable $\eta$ defined by 
$\eta=kt_1^{-1}\tau_f(t_1)),\, t_1=t/k,\, t\rightarrow 0$, $0<k<<t$ being an arbitrarily 
small random variable, so that O$(\eta^2)$=0, is a small scale microscopic variable when 
$\tau_f$ stands for the reflection symmetric fluctuating solution (\ref{fluc}).

\par  To summarize, we have presented two classes of nontrivial, exact solutions 
(\ref{ns}) and (\ref{fluc}), besides the standard solution (\ref{stand}), of the scale 
free equation (\ref{ss}), in the extended {\bf R}. The solution (\ref{fluc}) is 
infinitely differentiable and reflection symmetric, as is eq(\ref{stand}), but 
eq(\ref{ns}) is second derivative discontinuous and breaks the reflection symmetry 
spontaneously. Accordingly, the generalized solution (\ref{gen}) of eq(\ref{si}) is 
either reflection symmetric or otherwise, depending on the nature of the scale free 
component involving $\tau$ in its definition.      

\par Consequently, every real number is associated with a nontrivial neighbourhood of 
infinitesimally close elements ${\bf t}=t(1 + \eta)$, which would have been in a state 
of `equilibrium' random fluctuations, according to the definition 3, had there been only 
reversible, infinitely differentiable solutions (\ref{stand}) and (\ref{fluc}). Because 
of the second derivative discontinuous, asymmetric solution (\ref{ns}), the equilibrium, 
scale free fluctuation would be generically broken, inducing an irreversible, 
spontaneous evolution, albeit infinitely slow, to any approximate evaluation of the real 
number $t$. Note that any element in the set {\bf t} is an approximation of $t$, for a 
nonzero $\eta$. As shown in detail in Ref\cite{dp1,dp2}, this slow, intrinsic evolution 
would ultimately lead $\eta$ to approach the golden mean number ${\frac{\sqrt{5} 
-1}{2}}$, following the approximants of the corresponding continued fraction. ( The 
origin of the golden mean could be seen as follows. Being intrinsically an increasing 
variable (see below), the initially small variable $\eta$ would grow to a $\eta_- 
\lessapprox 1$. Following the solution (\ref{ns}), $\eta_-$ would then be replaced by 
$\eta_+^{-1}=1+\eta_1$,  $\eta_1$ being another small increasing variable and so on. The 
intrinsic evolution thus drives $\eta$ to follow the cascades of the golden mean 
continued fraction $1/(1+1/1+\ldots)$ ). The existence of the second derivative 
discontinuous solutions , with nontrivial dynamical properties, should be of crucial 
importance for future studies of dynamics. Our studies may only be considered as a few 
initial steps towards this direction.    

\par One may offer a `minimality argument' in support of the  inevitability of the time  
asymmetric solution in the extended solution space of eq(\ref{ss}). Note that in solving 
a first order (linear) differential equation (with $C^\infty$ data), the minimum 
condition one needs to meet is that the solution should at least be first order 
continuously differentiable. The infinitely differentiable solution (\ref{stand}) (and 
also (\ref{fluc})), however, corresponds to a situation where the information 
(knowledge) regarding the infinite differentiability of the solution is simply wasted. 
On the other hand, the time asymmetric solution just satisfies the minimal condition of 
being only the first derivative continuous. Now, invoking the common sense principle, 
viz., {\em the minimum information wastage (loss) means the maximum opportunity of 
progress}, one may assert that in the `universe' of all possible solutions of 
eq(\ref{si}), the time asymmetric solution would be naturally selected, because the said 
solution possesses the potentiality of growing into a more complex fractal like 
structure over a period of `infinitely' long time. Accordingly, the time asymmetric 
solution offers a dynamical system the maximum opportunity of evolution and growth.  We 
recall that the time sense is also generated here intrinsically by the solution itself. 
We note also that because of the intrinsic randomness and approximate nature of {\bf 
R}(=$R$), even the symmetric (fluctuating) solution (\ref{fluc}) fails to return an 
initial value $t_-$ (say) exactly, viz., $t_-\rightarrow t_+\rightarrow t_{-}^{\prime}$, 
$t_{-}^{\prime}$ being only approximately equal to $t_-$. This slight mismatch between 
the initial and final values would definitely induce a symmetry breaking, thereby 
realizing the time asymmetric solution in due course. Consequently, {\em the time 
asymmetric solution is uniquely realized in the extended scale free calculus}.   

\par Before closing this section, let us note that the above analysis reveals {\em an 
intrinsic mechanism} of change and evolution, which should universally be present not 
only in any dynamical system, but even rooted indelibly to the most fundamental of all 
mathematical entities, viz., the real number system . In other words, {\em the concepts 
of intrinsic evolution and time sense turn out to be truly   fundamental at par with  
real numbers, sets etc.} We note that in the ordinary framework of dynamical theories, 
motion and changes are meaningful only under the influence of an externally imposed 
force (or force field).  Consequently, the framework of the established dynamical 
theories such as Newtonian mechanics, Einstein's general relativity, quantum mechanics 
and so on, needs to be extended suitably incorporating this universally present 
intrinsic evolution \cite{egold,enash3,enash4}. More fundamentally, the study of the 
scale free calculus, as initiated in this series of works, may likely to provide a truly 
unified framework to understand natural and biological processes at different length, 
time and energy scales, in a scale invariant manner. 

\section{Dynamical numbers and cell division}

\par Natural numbers are most primitive of all numbers, being invented by primitive men 
to keep records of the number of objects they possessed. As such, natural numbers, and 
for that matter any number, are abstract mathematical concepts, which reside permanently 
and unchangeably in the `timeless' world of concepts. At an empirical level, however, a 
natural number $n$ (say) is associated with the cardinal number of the set of $n$ number 
of natural objects (we disregard the subtle distinction between the cardinal and ordinal 
numbers). For example, 2 is associated with the set of 2 billiard balls and so on. 
Instead of the set of billiard balls, one could have related 2 to the set of 2 flowers 
as well. Empirically, a set of living objects, e.g., a set of flowers, are prone to 
undergo rapid changes, relative to a set of inanimate objects, i.e., a set of billiard 
balls, for instance. Although, one tends habitually to disregard such natural mutations 
at the level of mathematical (arithmetical) analyses, the rigid, unchangeable structure 
of natural numbers are more closely respected by the sets of rigid, macroscopic material 
objects such as a set of billiard balls. Consequently, the traditional concept of 
natural numbers is lifeless, being devoid of an intrinsic potentiality of an 
irreversible change. Let us recall here that the theory of fuzzy sets and fuzzy logic is   
an approach to capture this `softer' (non-rigid) aspect of the natural systems. Our new 
results, even in the framework of the standard classical analysis offers an independent 
new approach to this problem.  

\par In view of the existence of intrinsically dynamic infinitesimals, every real 
number, and hence every natural number is now identified with a cell like structure (the 
intrinsically evolving set of infinitesimally close elements, any particular member of 
which being an approximate evaluation of the given natural number), which experiences a 
spontaneous fluctuation and irreversible evolution at every finer scale. The scale free 
number $1_f$, for example, now may be more accurately represented, at any given instance 
(as recorded by an ordinary clock), by the set of a new born living cell, $1_f=\{1 \, 
{\rm cell}\}$ (of course, $1_f$ can represent more complex sets such as a budding 
flower, or a sapling etc. We are, however, considering here the simplest case.). Now, 
$1_f=1 \pm \eta$, by definition, where $\eta=k t,\, k$ is an infinitesimal scale factor 
and $t$ is an O(1) intrinsic time variable which keeps record of the cell evolution. We 
note that $\eta$ satisfies eq(\ref{sc}), by definition, so that $k$ must solve 
eq(\ref{ss}) {\em in} $t^{-1}$. Had the fluctuations $\eta$ been reflection symmetric 
(i.e., $k$ be the reflection symmetric solution (\ref{fluc}), the cell like structure 
would have survived eternally. Because of the time asymmetric solution eq(\ref{ns}), 
$1_f$, however,  enjoys a number of evolutionary patterns, depending on the random 
selection of the $\pm$ sign, as $t\rightarrow 1/k$. 

\par Before detailing different evolutionary features, let us first examine more closely 
how the time asymmetric solution renders $\eta$ {\em dynamic}, with a slow, but directed 
evolution. As indicated, an infinitesimal $\eta$ can always be parametrized as 
$\eta=kt$, where $t=1 \pm \eta _1,\, \eta_1 >0$ being another infinitesimal. As seen 
already, $\eta _1$ (and hence $\eta$) cannot be static, ordinary constant, and 
consequently, $t$ must undergo `infinitely' slow, irreversible evolution towards $1/k$, 
because of the cooperative effect of the infintesimally small scale free (random) 
fluctuations, represented by the time asymmetric solutions (\ref{ns}) and (\ref{as}). In 
fact, the solution (\ref{as}) forces $t$ to change irreversibly  towards the forward 
direction by the applications of infinitesimal scalings, inversions and translations. 
Let for instance, $t=1-\eta_1$. Then infinitesimal fluctuations transform $t$ to 
$t_1=t^{-\alpha}=1+\alpha\eta_1 \equiv 1+\eta_2 ,\,\alpha>0$ (by a combination of 
inversion and scaling). This new value of $t_1$ could now continue to grow by analogous 
applications of scalings of the form $t_1 \rightarrow t_2=t_1^{\alpha_1}$ and so on, 
until the successively transformed value of $\eta_n$ after $n$th iteration becomes of 
the order of a small macroscopic variable $\bar t (0<\bar t<< 1)$, which will then 
increase monotonically (and linearly) till $t=1+\bar t \sim 1/k $. Analogously, 
$t=1+\eta_1$ would also increase towards $1/k$ by successive applications of 
infinitesimal scalings and translations. This intrinsic monotonic movement of $t$ must 
be infinitely slow, since the number of iterations $n$ should be an infinitely large 
positive integer, for any finite $\alpha$.

\par Once the intrinsic time like variable $t $ grows to O(1/k), so that $\eta=1 \pm 
\eta^{\prime}$, the cell like structure of of $1_f$ would either (a) spontaneously 
replicate to two daughter cells: $1_f (=1+\eta) \rightarrow 2_f= 1_f^{\prime} + 
1_f^{\prime}$, where $1_f^{\prime}=1+{\frac{1}{2}\eta}^{\prime}$ is an almost exact 
replica of $1_f$, but for an infinitesimal variation, or (b) annihilate $1_f (=1-\eta) 
\rightarrow 0_f$. The second possibility may signify the `death' of the new born cell, 
whereas the first case corresponds to the replication of the mother cell into daughter 
cells by the process of mitotic division. Because of the random occurrence of the $\pm$ 
sign, the process of celluar  division would likely to be halted after a finite number 
of cell replication when the growing mode (viz., $1_f$ with a `+' sign) meets for the 
first time the `-' sign and is thus transferred to the decaying mode. (The role of the 
golden mean, and the Fibonacci numbers in  cell division would be examined separately. ) 

\par To summarize, the scale free natural numbers seems to enjoy living cell like 
intrinsic dynamical behaviours. One needs, however, to be reminded that the present 
analysis cannot pretend to be a complete description of the cell division process. The 
significance of the present discussion is to highlight the importance of the {\em 
intrinsic dynamics} in the contemporary studies of cell division and other biological 
processes. In fact, one may call a system living whose dynamics is more efficiently 
controlled by the scale free intrinsic time like component, compared to mere external 
time treatment, as in vague in the contemporary studies. In this sense, every system is 
living, in the present extended framework of scale free dynamics. However, one system 
may be more living than other. For example, a billiard ball is material (nonliving)    
since the intrinsic time and corresponding intrinsic dynamics become relevant only 
infinitely slowly, making it meaningless in the context of ordinary Newtonian dynamics. 
However, for a living biological cell, the intrinsic dynamics should  become effective 
even at an order O(1) of the ordinary time scale. 

\section{Newton's first law of motion: An intelligent version}

\par The first law of motion in classical mechanics states that ` a particle in uniform 
rectilinear motion or in rest will continue to remain in the state of uniform motion or 
rest unless disturbed by an {\em externally applied} force'. A direct consequence of 
this law is the following. Let A and B be two points, moving uniformly with speed 1 unit 
along $x$ axis from the initial points $x=0$ and $x=2$ respectively, in the opposite 
directions. Then, according to the first law, the two particles would definitely collide 
at $x=1$.

\par However, an interesting variation of the first law could be possible in the 
extended scale free mechanics, incorporating time asymmetric solutions. We note at first 
that the instantaneous positions of the two particles are given by $x_A=t$ and $x_B=2-t$ 
respectively, $t$ being the time variable. Thus as $t \rightarrow 1^-$, $x_A=1-\eta$ and 
$x_B=1+\eta$, $0<\eta<< 1$. Consequently, for an infinitesimal $\eta$, A and B may avoid 
direct collision, by simply flipping the sign of $\eta$, thereby interchanging the role 
of $x_A$ and $x_B$. As a result, A and B, as they approach towards $x=1$, would jump 
around each other, and then continue their respective motions uninterrupted, viz., A 
towards $x=2$ and B towards $x=0$. One is, therefore, tempted to infer that the time 
asymmetric solution  instills, so to speak,  {\em intelligence} into moving particles, 
leading to an avoidance of the otherwise unavoidable collision. Further, the reason of 
this avoidance is not due to some external forces, but solely an {\em intrinsic} 
property of the asymmetric solutions of the linear differential equations. The avoidance 
of an immanent  collision, as dictated by the Newton's first law, provides an instance 
of a `controlled' behaviour, which along with the intrinsic irreversibility, as pointed 
out in Sec.2.2, may be considered as the hallmark of an intelligent living system.     

\section*{Appendix}

\par An exact (to all orders of $\eta$ ) random solution of eq(\ref{ss}) is obtained 
when the  approximate solution (\ref{ns}) is improved recursively by self-similar 
correction factors over smaller and smaller scales $\eta,\, \eta^2,\,\eta^4, \, \ldots$.
Let 

$${
\tau_f(t)=
\cases{
\tau_- & {\rm if} $t\lessapprox 1$\cr
\tau_+ & {\rm if} $t\gtrapprox 1$},\,\,
\tau_{-}(t_-)=(1/t_+)f_-(\eta),\, \tau_+(t_+)=t_+
}\eqno{({\rm A}1)}$$
 
 \noindent We consider the nontrivial part only. One obtains

$${
t_-{\frac{{\rm d}\tau_{-}}{{\rm d} t_-}}=\tau_- ({\frac{t_-}{t_+}}-t_-
{\frac{f_-^\prime}{f_-}})
}\eqno{({\rm A}2)}$$

\noindent where $f^\prime={\frac{{\rm d}f}{{\rm d}\eta}}$. As a result, $\tau_-$ solves 
eq(\ref{ss}) exactly provided $f_-$ solves exactly the self-similar equation   

$${
t_-^\prime{\frac{{\rm d}f_-}{{\rm d} t_-^\prime}}=f_- 
}\eqno{({\rm A}3)}$$

\noindent in the smaller scale variable $\ln(1/ t_-^\prime)= \ln (1-\eta^2)^{-1}$. The 
exact solution is thus obtained as an infinite product 

$${
\tau_={\frac{1}{t_+}}{\frac{1}{t_+^\prime}}{\frac{1}{t_+^{\prime\prime}}}\ldots
}\eqno{({\rm A}4)}$$

\noindent where $ t_+^\prime=  (1+\eta^2),\, t_+^{\prime\prime}= (1+\eta^4)$ etc. One 
therefore reproduces the standard solution $\tau_-=1-\eta$ (since $t_+t_+^\prime\ldots 
=(1-\eta)^{-1}$), when $\eta$ is an ordinary real variable with exact values. In the 
present extended treatment $\eta$, however, is a small random variable. The solution 
$\tau_-$ could therefore be interpreted as an infinite multiplicative process, the 
self-similar factors of which may be treated as independent random variables defined 
recursively via eq(\ref{ss}) and eq({\rm A}4). Maintaining continuity of the second 
derivative of $\tau_f$ at $\eta=0$, in the ordinary real calculus sense, would be 
impossible.

\end{document}